\documentclass[daat]{ipart}

\RequirePackage{hyperref}
\usepackage{graphicx}
\usepackage{subcaption}
\usepackage{algorithm}
\usepackage{algorithmic}
\usepackage{booktabs}

\startlocaldefs
\theoremstyle{plain}

\newtheorem{remark}{Remark}

\endlocaldefs
\startlocaldefs
\theoremstyle{plain}

\endlocaldefs

\pubyear{2025}
\volume{1}
\issue{1}
\firstpage{1}
\lastpage{4}

\begin{document}

\begin{frontmatter}
\title{A novel implementation of Yau-Yau filter for time-variant 
nonlinear problems}

\begin{aug}
    \author{\inits{Y.}\fnms{Yuzhong} \snm{Hu}\ead[label=e1]{huyuzhong@bimsa.cn}},
    \address{Beijing Key Laboratory of Topological Statistics and Applications for Complex Systems, Beijing Institute of Mathematical Sciences and Applications, \\
             Beijing 101408, China\\
             Chinese Academy of Mathematics and System Science,\\
             Beijing 100190, China\\
             \printead{e1}}
    \author{\inits{J.}\fnms{Jiayi} \snm{Kang}\ead[label=e2]{kangjiayi@bimsa.cn}},
    \address{Beijing Key Laboratory of Topological Statistics and Applications for Complex Systems, Beijing Institute of Mathematical Sciences and Applications, \\
             Beijing 101408, China\\
             \printead{e2}}
    \author{\inits{L.}\fnms{Lei} \snm{Ma}\ead[label=e3]{malei@bimsa.cn}},
    \address{Beijing Key Laboratory of Topological Statistics and Applications for Complex Systems, Beijing Institute of Mathematical Sciences and Applications, \\
             Beijing 101408, China\\
             Chinese Academy of Mathematics and System Science,\\
             Beijing 100190, China\\
             \printead{e3}}
	\and
    \author{\inits{X.}\fnms{Xiaoming} \snm{Zhang}\thanksref{t2}\ead[label=e4]{zhangxiaoming@bimsa.cn}}
    \address{Beijing Key Laboratory of Topological Statistics and Applications for Complex Systems, Beijing Institute of Mathematical Sciences and Applications, \\
             Beijing 101408, China\\
             \printead{e4}}
    \thankstext{t2}{Corresponding author.}
\end{aug}
%

\begin{abstract}
Nonlinear filter has long been an important problem in practical industrial applications. The Yau-Yau method is a highly versatile framework that transforms nonlinear filtering problems into initial-value problems governed by the Forward Kolmogorov Equation (FKE). Previous researches have shown that the method can be applied to highly nonlinear and high dimensional problems. However, when time-varying coefficients are involved in the system models, developing an implementation of the method with high computational speed and low data storage still presents a challenge. To address these limitations, this paper proposes a novel numerical algorithm that incorporates physics-informed neural network (PINN) and principal component analysis (PCA) to solve the FKE approximately. Equipped with this algorithm, the Yau-Yau filter can be implemented by an offline stage for the training of a solver for the approximate solution of FKE and an online stage for its execution. Results of three examples indicate that this implementation is accurate, both time-efficient and storage-efficient for online computation, and is superior than existing nonlinear filtering methods such as extended Kalman filter and particle filter. It is capable of applications to practical nonlinear time-variant filtering problems. 
\end{abstract}



\end{frontmatter}

\section{Introduction}
Nonlinear filtering has long been a much studied discipline within the broader field of signal processing and control theory, due to its importance in many applications such as robotics control, target tracking and maneuvering, electrical power systems, financial market prediction and management, biomedical signal processing, and vehicle navigation \cite{Chen2009, Date2011, Zhang2011}. In contrast to linear filtering where the system equations are linear, nonlinear filtering addresses more complex problems where the system dynamics models and observation models are inherently nonlinear. The main challenge in nonlinear filtering lies in estimating the system’s state over time from a sequence of non-Gaussian distributed noisy observations, where both the state transition and observation models are nonlinear \cite{Crisan2011}.

The origins of filtering theory date back to Gauss's early work almost two centuries ago, followed by Wiener’s substantial contributions \cite{Kutschireiter2020}. In 1960, Rudolf E. Kalman introduced the Kalman filter \cite{Kalman1960, Kalman1961}, an optimal filtering technique for linear systems with Gaussian noise. However, the direct application of the Kalman filter is not suitable for nonlinear systems. Over time, several extensions have been developed to address these limitations, including the extended Kalman filter (EKF) \cite{Jazwinski2007}, unscented Kalman filter (UKF) \cite{Julier1997}, ensemble Kalman filter (EnKF) \cite{Evensen2003}, and particle filters (PF) \cite{Gordon1993}. EKF linearizes nonlinear dynamics via first-order Taylor expansions around current estimates, providing efficiency but potentially significant errors in strongly nonlinear contexts. The UKF utilizes unscented transformations to achieve higher accuracy at the cost of greater computational complexity. EnKF handles nonlinearity efficiently through an ensemble of state estimates without full covariance matrix calculations, making it suitable for high-dimensional problems such as weather forecasting. Meanwhile, Particle Filters represent posterior distributions with weighted particles, offering flexibility for strong nonlinearities but suffering from computational intensity due to sample degeneracy.

Besides Kalman-type filters and PF, it is also possible to solve the filtering problem by tracking the evolution of the probability density function (PDF) of the state process given the observations. In the 1960s, Duncan \cite{Duncan1967}, Mortensen\cite{Mortensen1966}, and Zakai\cite{Zakai1969} independently derived the so-called Duncan–Mortensen–Zakai’s (DMZ) equation, a stochastic partial differential equation (SPDE), which is satisfied by the unnormalized conditional density function of the states. In general filtering systems, the DMZ equation lacks a closed-form solution, and solving it for more complex systems has become a central challenge in nonlinear filtering theory. In 2008, Yau and Yau \cite{Yau2000,Yau2008} made significant advancements by solving the DMZ equation under the least relaxed conditions, introducing a two-stage method for its solution in which an on-line computation consists only of an explicit exponential transformation which is time-efficient, and an offline procedure for the expensive computation for solving the Forward Kolmogorov Equation (FKE). They also provide a theoretical guarantee of algorithmic convergence.

The Yau-Yau method successfully transforms the challenges of solving filtering problems into computational mathematical problems, specifically the Cauchy problem of solving FKE. Inspired by subsequent methods in computational mathematics, implementations such as the Hermite spectral method \cite{Luo2013}, Legendre-Galerkin spectral method \cite{Dong2020}, Euler-optimization methods \cite{Shi2018}, and quasi-implicit Euler techniques \cite{Yueh2014} have been proposed to solve FKE offline. Nevertheless, traditional grid-based numerical methods like finite difference and finite element approaches require extensive computational resources and face the "curse of dimensionality" when applied to high-dimensional problems. Spectral methods similarly encounter dimensionality challenges. To address these computational burdens, reduced-order methods such as proper orthogonal decomposition and deep Galerkin approaches have been introduced \cite{Wang2019,Shi2024}. 

In recent years, the advent of machine learning techniques, particularly neural networks, has been developed for solving differential equations and thus opens new avenues for addressing the nonlinear filtering problem. Chen et al \cite{Chen2024} developed a computationally efficient method for nonlinear filter design and recurrent neural networks. Tao et al \cite{Tao2023} proposed a neural projection filter to handle high dimensional nonlinear filter problems with robustness and efficiency. Physics-Informed Neural Networks (PINNs) \cite{Raissi2019} have obtained significant attention for their capability to incorporate physical laws such as in the form of differential equations directly into the learning process, ensuring that the solutions adhere to the underlying governing system dynamics. PINN’s theoretical basis and generalized error bounds estimation have been studied \cite{Mishra2023,Chen1995}. It offers a promising alternative to mesh-based and spectrum-based methods, particularly in scenarios when dealing with highly nonlinear and high dimensional systems. 

However, existing research primarily focusses on time-invariant systems, where the system model and the observation model do not have time-variant coefficients. Yet, many practical industrial problems involve explicit time dependency in the state process and the observation process \cite{Speidel2021, Pintelon2016, Olaleye2004}, introducing time-variant coefficients in the FKE. Solving FKE under such time-variant conditions necessitates substantial computational resources and memory for storing pre-calculated solutions. Therefore, a significant research gap remains: developing accurate, computationally efficient, and storage-efficient methods for high-dimensional, time-variant nonlinear filtering systems.

Direct application of PINN to time-variant filtering problems is still a challenge since FKE involves time-dependent coefficients whose values are different for different time intervals. It is not feasible to store all the FKE solutions for all kinds of coefficients in the storage for online use. To address this challenge, dimensional reduction methods such as principal component analysis (PCA) \cite{Abdi2010} can be a candidate to reduce the information we need to store. 

In this paper, we propose a novel method that integrates Yau-Yau filter, PINN, and PCA to produce real-time predictions of states for nonlinear filtering problems. We refer to this new method as the PINN-based Yau-Yau filter (PINNYYF). The implementation comprises two stages: an offline stage (Stage I) subdivided into Stage IA—solving the FKE with PINNs given initial conditions—and Stage IB—constructing an approximate FKE solver by combining PINNs with PCA; and an online stage (Stage II), which performs near-real-time predictions of the state variables through the utilization of the approximate FKE solver.

The main contributions of this paper are listed as follows:
\begin{enumerate}
\item We apply PINN for the solution of FKE in the Yau-Yau filter framework for time-variant systems. The numerical results demonstrate that the application of PINN can achieve high accuracy.
\item Development of an interpretable FKE solver using PINNs and PCA, trained with initial condition-solution pairs at each time step. This solver significantly improves computational efficiency and storage efficiency in time-variant nonlinear filtering.
\end{enumerate}

This paper is structured as follows: Section 2 briefly describes the nonlinear filtering problem, the Yau-Yau method and the Physics Informed Neural Network. Section 3 proposes the two-stage implementation of the Yau-Yau method for time-variant systems. In Section 4, we apply PINNYYF to three examples illustrating the application procedure and demonstrating its effectiveness in comparison with EKF and PF. Section 5 concludes the paper with a discussion on the implications of our findings and potential directions for future research.

\section{Preliminaries}

\subsection{Nonlinear Filtering Problem and the Yau-Yau Filter}

A signal-observation model can be represented by a system of stochastic differential equations:

\begin{equation}
\begin{aligned}
d\mathbf{x}_t &= f(\mathbf{x}_t,t)dt + G(\mathbf{x}_t,t)d\mathbf{w}_t, \\
d\mathbf{y}_t &= h(\mathbf{x}_t,t)dt + d\mathbf{v}_t,
\end{aligned}
\end{equation}
where \( \mathbf{x}(t) \) is the d-dimensional vector of the state at time \( t \), and \( \mathbf{y}(t) \) is the m-dimensional vector of the observation at time \( t \). Here, \( f \) and \( h \) are vector-valued functions, \( G \) is a $d\times r$ matrix-valued function, and $v_t$ is a r-dimensional Brownian motion process with $\mathbb{E}(\mathbf{v}_t\, \mathbf{v}_t^\top) = Q(t)\,dt, \quad Q(t) \in \mathbb{R}^{r \times r}.$ The processes \( \mathbf{w}(t) \) and \( \mathbf{z}(t) \) are standard Brownian motions.Besides, $y_t$ and $h$ are $m$-vectors, and $w_t$ is an $m$-dimensional Brownian motion process with
$
\mathbb{E}(d\mathbf{w}_t\, d\mathbf{w}_t^\top) = S(t)\,dt, \quad S(t) \in \mathbb{R}^{m \times m}.
$
The initial state $\mathbf{x}_0$ follows an initial distribution $\sigma_0(x)$, and we assume $\mathbf{v}_t$, $\mathbf{w}_t$, and $\mathbf{x}_0$ are mutually independent.

Let $F_t = \sigma\{y_s : 0 \leq s \leq t\}$ denote the smallest $\sigma$-field generated by historical observations. Denote $p(x,t\,|\,F_t)$ as the conditional probability density of state $x_t$, and in the following context denote it as $p(x,t)$ for short. $p(x,t)$ can be obtained through the unconditional probability density $\sigma(x,t)$ by solving the DMZ equation \cite{Duncan1967,Mortensen1966,Zakai1969} given by:

\begin{equation}
\begin{cases}
d\sigma(x,t) = \mathcal{L} \sigma(x,t)\,dt + \sigma(x,t)\, h^\top(x)\, S^{-1} \,dy(t), \\
\sigma(x,0) = \sigma_0(x)
\end{cases}
\label{eq:DMZ}
\end{equation}

where
\begin{equation}
\mathcal{L}(*) := \frac{1}{2} \sum_{i,j=1}^n \frac{\partial^2}{\partial x_i \partial x_j} \left[ (GQG^\top)_{ij}\, * \right] - \sum_{i=1}^n \frac{\partial}{\partial x_i} \left( f_i\, * \right),
\end{equation}

and $p(x,t)$ can be obtained by
$p(x,t) = \frac{\sigma(x,t)}{\int \sigma(x,t)\,dx}$.
Firstly, we apply an invertible exponential transformation \cite{Davis1980} for each given observation,
\begin{equation}
\sigma(x,t) = \exp\left[ h^\top(x,t)\, S^{-1} y_t \right] \rho(x,t).
\end{equation}

Equation \ref{eq:DMZ} can be transformed into a deterministic partial differential equation with stochastic coefficients, which is referred to as the “pathwise-robust” DMZ equation:

\begin{equation}
\label{eq:robust}
\begin{cases}
\frac{\partial \rho(x,t)}{\partial t}
+ \frac{\partial}{\partial t} \left\{ (h^\top S^{-1})^\top y_t \right\} \rho(x,t)
\\= \exp\left\{ -h^\top S^{-1} y_t \right\}
\left( L - \frac{1}{2} h^\top S^{-1} h \right) \\
\cdot \exp\left\{ h^\top S^{-1} y_t \right\} \rho(x,t), \\
\rho(x,0) = \sigma_0(x)
\end{cases}
\end{equation}

The observations arrive at discrete times $\tau_i = i\Delta t$, $i = 0, 1, 2, \ldots, N_T$, and $\Delta t = T/N_T$. Denote $\rho_i$ as the solution of the robust DMZ equation (Eq.~\ref{eq:robust} with $y_t = y_{\tau_{i-1}}$ on the time interval $t \in [\tau_{i-1}, \tau_i]$, $i = 1, 2, \ldots, N_T$, i.e.,

\begin{equation}\label{eq:transformed}
\begin{cases}
\frac{\partial \rho_i}{\partial t}(x, t) + \frac{\partial}{\partial t}\{(h^TS^{-1})^T y_{\tau_{i-1}}\}\rho_i(x, t) \\
= \exp\{-h^T S^{-1} y_{\tau_{i-1}}\}\left(L - \frac{1}{2}h^T S^{-1} h\right)\cdot \\
\exp\{h^T S^{-1}y_{\tau_{i-1}}\}\rho_i(x, t) \\
\rho_1(x, 0) = \sigma_0(x) \\
\rho_i(x, \tau_{i-1}) = \rho_{i-1}(x, \tau_{i-1}) \text{ for } i = 2, 3,\ldots,N_T.
\end{cases}
\end{equation}

Putting together the $\rho_i$'s yields an approximation to $\rho$ in Eq.~\ref{eq:robust}, i.e.,
\begin{equation}
\rho(x,t) \approx \sum_{i=1}^{N_T} \chi_{[\tau_{i-1}, \tau_i]}(t)\, \rho_i(x,t),
\end{equation}
where $\chi_{[a,b]}(t)$ is the indicator function of the interval $[a,b]$, i.e., $\chi_{[a,b]}(t) = 1$ if $t \in [a,b]$ and $0$ otherwise.

With an inverse transformation, $u_i(x,t) = \exp\left[ h^\top\, S^{-1}\, y_{\tau_{i-1}} \right] \rho_i(x,t)$, Eq.~(6) satisfies the Forward Kolmogorov Equation (FKE) \cite{Yau2008}:
\begin{equation}\label{eq:FKE}
\frac{\partial u_i}{\partial t}(x,t) = \left( \mathcal{L} - \frac{1}{2} h^\top S^{-1} h \right) u_i(x,t),
\end{equation}
where $\mathcal{L}$ is defined in Eq.~(3), and the initial conditions are as follows:
\begin{equation}
\begin{cases}
u_1(x,0) &= \sigma_0(x), \text{or}\\
u_i(x, \tau_{i-1}) &= 
\exp\left[ h^\top\, S^{-1}\, \left( y_{\tau_{i-1}} - y_{\tau_{i-2}} \right) \right]\\ &\cdot u_{i-1}(x, \tau_{i-1}), \quad i \geq 2.
\end{cases}
\end{equation}
The FKE in Eq.~\ref{eq:FKE} is obviously independent of observations and can therefore be solved offline given initial conditions to reduce the online computational burden. The implementation procedure of the Yau-Yau method consists of two stages:

\textbf{Stage I (Offline)}: $u_i(x,\tau_{i-1}) \rightarrow u_i(x,\tau_i)$

Solve the FKE (Eq.~(8)) over the time interval $[\tau_{i-1}, \tau_i]$ with initial value $u_i(x,\tau_{i-1})$. This stage can be computed offline for any given initial values since Eq.~(8) is independent of the observations.

\textbf{Stage II (Online)}: $u_i(x,\tau_i) \rightarrow u_{i+1}(x,\tau_i)$

When the new observation $y_{\tau_i}$ arrives at time $\tau_i$, the initial value of $u_{i+1}(x,t)$ for the next time interval $[\tau_i, \tau_{i+1}]$ is updated. For $t \in [0, \tau_1]$, the initial condition is
\[
u_1(x,0) = \sigma_0(x).
\]
For $i \geq 2$, the update of the initial value is given by:
\begin{equation}
u_i(x, \tau_{i-1}) = \exp\left[ h^\top\, S^{-1}\, \left( y_{\tau_{i-1}} - y_{\tau_{i-2}} \right) \right] \cdot u_{i-1}(x, \tau_{i-1}).
\end{equation}

With these two computational stages for each of all time intervals, the functions $u_{i+1}(x, \tau_i)$ serve as good approximations to the solution $\sigma(x, \tau_i)$ of the original DMZ equation (Eq.~(2)) at time $t = \tau_i$.

\subsection{Physics-Informed Neural Networks}

PINNs are an innovative class of machine learning models that incorporate physical laws, typically expressed as partial differential equations (PDEs), directly into the training process of neural networks. Unlike traditional neural networks, which rely solely on data to learn patterns and make predictions, PINNs leverage both data and known physical principles in terms of differential equations, ensuring that the model’s outputs adhere to the underlying physics of the problem.

\begin{figure}[htbp]
    \centering
    \includegraphics[width=\linewidth]{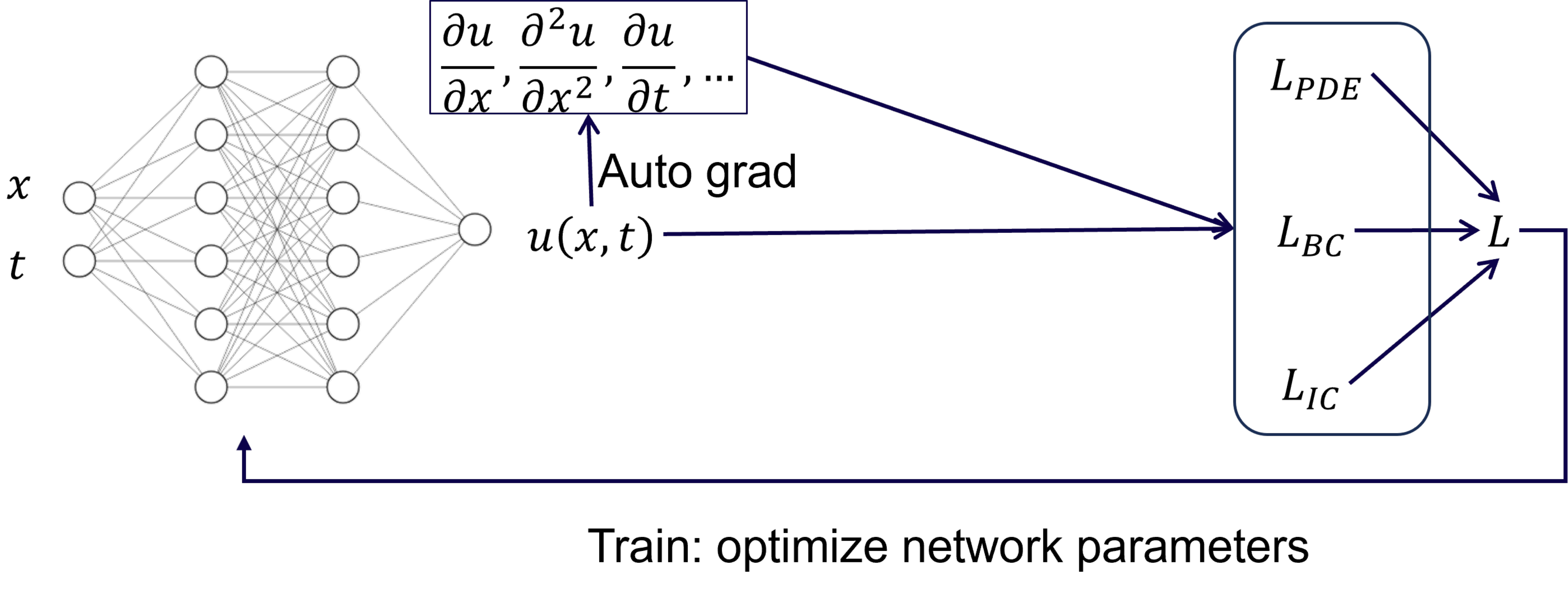} 
    \caption{The general structure of a Physics-Informed Neural Network (PINN).}
    \label{fig:pinn_structure}
\end{figure}

Figure~\ref{fig:pinn_structure} illustrates the general structure of a PINN. The network takes as input the spatiotemporal coordinates $(x,t)$, referred to as collocation points, and computes an approximation $u(x,t)$ of the true solution to the target problem. The outputs of the network, along with their required derivatives, are used to construct the loss function. These derivatives are efficiently evaluated using automatic differentiation techniques \cite{Paszke2017}.

The loss function typically consists of multiple components. For example, $L_{\text{PDE}}$ measures the residuals of the governing partial differential equation, $L_{\text{IC}}$ enforces the initial conditions, and $L_{\text{BC}}$ corresponds to the boundary conditions. The total loss is formed by a weighted sum of these terms:
\begin{equation}
L(\theta) = \lambda_{\text{PDE}} L_{\text{PDE}} + \lambda_{\text{IC}} L_{\text{IC}} + \lambda_{\text{BC}} L_{\text{BC}},
\end{equation}
which is minimized iteratively to optimize the network parameters $\theta$.

As a PDE solver, PINN can handle high-dimensional problems with complex geometries without requiring mesh generation or discretization, which are typically necessary in traditional methods. In the solution of FKE of each time step, the network can be obtained by finetuning the network for the previous time step, and is thus computationally efficient. 

\subsection{Principal Component Analysis}

PCA \cite{Abdi2010} is widely used for dimensionality reduction and feature extraction in high-dimensional datasets. PCA transforms data into orthogonal principal components by computing the covariance matrix and identifying eigenvectors associated with its largest eigenvalues. These components capture maximum variance, significantly reducing dimensionality while preserving essential information. PCA projects the original dataset onto selected principal components, creating a low-dimensional representation suitable for efficient storage and rapid computations.

\section{Implementation Algorithm for Time-variant Systems}

Our algorithm addresses two key challenges inherent in nonlinear filtering for time-variant systems: First, constructing a computationally efficient and storage efficient FKE solver; and second, incorporating time-dependent features in the filtering process. 

Traditional numerical methods for solving the FKE such as finite difference method often require substantial computational resources and large storage capacities. Additionally, spectral methods necessitate a predetermined set of basis functions, whose selection critically influences the performance of the resulting filter. Determining the optimal choice of these basis functions remains a non-trivial and problem-specific task. Although data-driven filtering techniques, such as the RNNYYF \cite{Chen2024} and neural projection filters \cite{Tao2023}, have demonstrated significant progress in computational efficiency, these methods typically do not explicitly solve the underlying FKE. Furthermore, existing approaches based on the traditional Yau-Yau method often requires substantial storage or computational demands for modeling systems with explicit time-variant coefficients \cite{Dong2020, Chen2017}. 

To tackle these challenges, we introduce a two-stage framework:

\begin{itemize}
    \item \textbf{Stage I (offline):}
    \begin{itemize}
        \item Stage IA: Train a PINN-based FKE solver.
        \item Stage IB: Construct an approximate time-variant FKE solver using PCA.
    \end{itemize}
    \item \textbf{Stage II (online):} 
    Execute the approximate FKE solver and obtain the end time solutions.
\end{itemize}

\subsection{Stage IA (offline): Train FKE solver with PINN}

In this offline substage, we employ PINN to solve the initial value problem of the FKE. For each time step $[\tau_{i-1}, \tau_i]$ and any initial condition $u_i(x, \tau_{i-1})$, $x \in \Omega$, we train a PINN that takes the space-time coordinates $(x,t) \in \Omega \times [\tau_{i-1}, \tau_i]$ as input, and output $\hat{u}(x,t)$ that approximates the FKE solution $u_i(x,t)$. The training objective is formulated as a composite loss function consisting of three distinct components:
\begin{equation}
L = L_{\text{FKE}} + \lambda_{\text{IC}} L_{\text{IC}} + \lambda_{\text{BC}} L_{\text{BC}}
\end{equation}
where,
\begin{itemize}
    \item 
    $L_{\text{FKE}} = \frac{1}{N_{\text{FKE}}} \sum_{j=1}^{N_{\text{FKE}}} \text{FKE}(u_i(x^{(j)}_{\text{FKE}}, t^{(j)})) 
    = \frac{1}{N_{\text{FKE}}} \sum_{j=1}^{N_{\text{FKE}}} 
    \left[ \frac{\partial \hat{u}}{\partial t}- \left( \mathcal{L} - \frac{1}{2} h^\top S^{-1} h \right) \hat{u} \right]^2 \bigg|_{(x^{(j)}_{\text{FKE}}, t^{(j)})}$

    measures the residuals of FKE at $N_{\text{FKE}}$ collocation points $(x^{(j)}_{\text{FKE}}, t^{(j)})$, $j = 1, 2, \ldots, N_{\text{FKE}}$, randomly selected from $\Omega \times [\tau_{i-1}, \tau_i]$, to enforce the governing FKE.
    
    \item 
    $L_{\text{IC}} = \frac{1}{N_{\text{IC}}} \sum_{j=1}^{N_{\text{IC}}} \text{IC}(u_i(x^{(j)}_{\text{IC}}, \tau_{i-1})) 
    = \frac{1}{N_{\text{IC}}} \sum_{j=1}^{N_{\text{IC}}} 
    \left[ \hat{u} - u_i \right]^2 \bigg|_{(x^{(j)}_{\text{IC}}, \tau_{i-1})}$

    enforces the initial condition through $N_{\text{IC}}$ collocation points $x^{(j)}_{\text{IC}}$, $j = 1, 2, \ldots, N_{\text{IC}}$, randomly selected from $\Omega$.

    \item 
    $L_{\text{BC}} = \frac{1}{N_{\text{BC}}} \sum_{j=1}^{N_{\text{BC}}} \text{BC}(u_i(x^{(j)}_{\text{BC}}, t^{(j)})) 
    = \frac{1}{N_{\text{BC}}} \sum_{j=1}^{N_{\text{BC}}} 
    \left[ \hat{u} \right]^2 \bigg|_{(x^{(j)}_{\text{BC}}, t^{(j)})}$

    enforces zero boundary condition, evaluated over $N_{\text{BC}}$ collocation points $(x^{(j)}_{\text{BC}}, t^{(j)})$, $j = 1, 2, \ldots, N_{\text{BC}}$, randomly selected from $\partial \Omega \times [\tau_{i-1}, \tau_i]$.
\end{itemize}

Parameters $\lambda_{\text{IC}}$, $\lambda_{\text{BC}}$ are the weights of the loss terms, balancing the relative importances of these three losses, while $N_{\text{FKE}}, N_{\text{IC}}, N_{\text{BC}}$ are selected to ensure high accuracy.

Before training the network, we generate a set of the state and observations $x_t, y_t$ with Monte-Carlo simulations for $t \in [0, T]$. At each time step, the neural network is trained to obtain the end time solution \( u_i(x, \tau_i) \), and the corresponding simulated observations \( y_{\tau_i} \) are used to update the initial condition. This iterative process continues to generate pairs of initial conditions and solutions.

\begin{figure}[htbp]
    \centering
    \includegraphics[width=\linewidth]{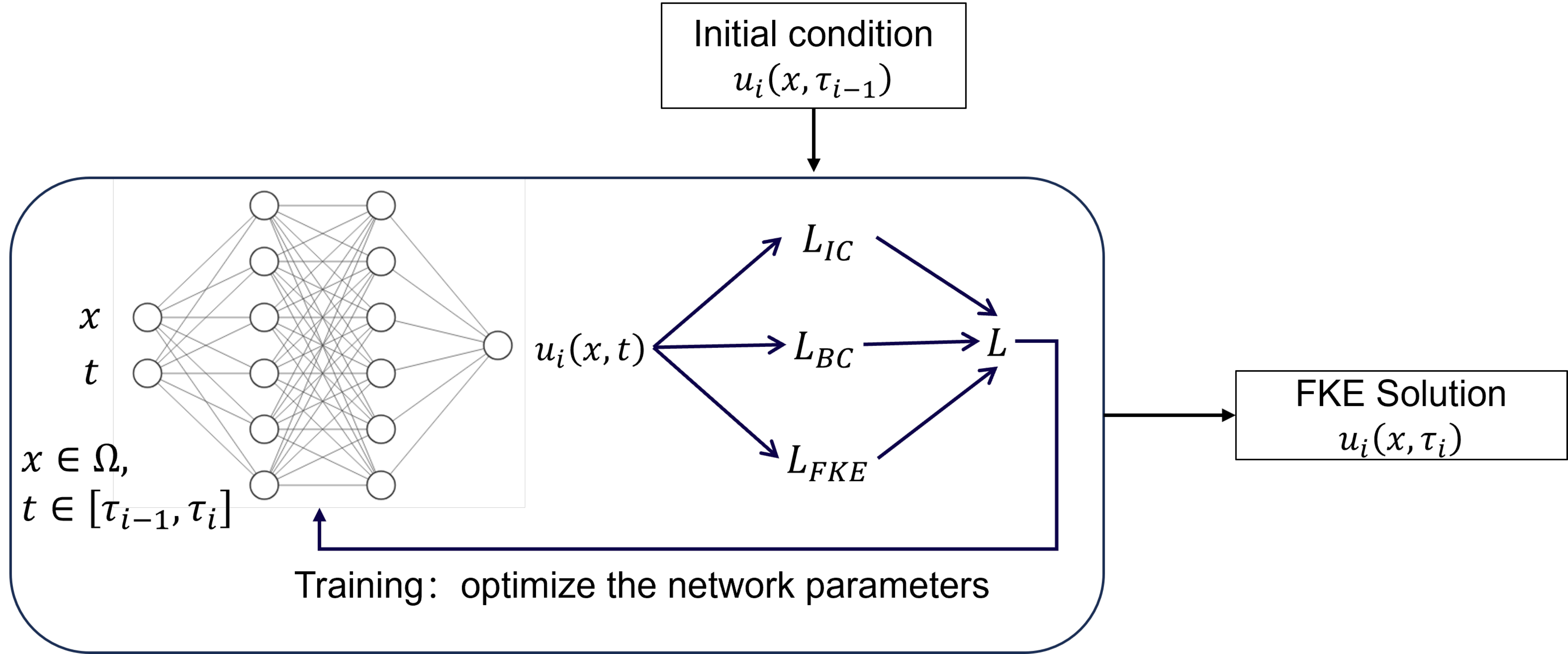} 
    \caption{The training process of FKE solver with PINN.}
    \label{fig:stageIA}
\end{figure}

 To accelerate training, we initialize the network parameters using weights from the preceding time step. For the initial time step $t \in [\tau_0, \tau_1]$, the network parameters are randomly initialized. Fig.~\ref{fig:stageIA} describes the general training process of the FKE solver.

The general training algorithm of PINNs corresponding to Stage IA is as follows:
\begin{algorithm}[htbp]
\caption{Train FKE solver with PINN (offline)}
\label{alg:stageIA}
\begin{algorithmic}[1]
\STATE \textbf{Initialization:} given the domain $\Omega$, time interval $\Delta t$, number of time steps $N_T$, convergence threshold $\epsilon$, and maximum number of epochs $E_{\max}$, the observation $y_{\tau_i}$ for $i = 0, 1, 2, \ldots, N_T$. We randomly initialize the network parameters $\theta$.
\FOR{$i = 1$ to $N_T$}
    \STATE Initialize the network parameters $\theta$ from preceding time step
    \FOR{epoch = 1 to $E_{\max}$}
        \STATE Compute the loss $L$ using Eq.~(10), and optimize the network parameters $\theta$
        \IF{$L < \epsilon$}
            \STATE break the inner for loop
        \ENDIF
    \ENDFOR
    \STATE Use the network to evaluate the terminal solution $u_i(x, \tau_i)$, and update the initial condition $u_{i+1}(x, \tau_i)$ with upcoming observations $y_{\tau_i}$
\ENDFOR
\end{algorithmic}
\end{algorithm}

\subsection{Stage IB (offline): Train the approximate time-variant FKE solver}

In Stage IA, we obtain the initial conditions $u_i(x,\tau_{i-1})$ and the corresponding end time solutions $u_i(x,\tau_i)$. We then use PCA to extract the top $K$ principal components from these functions, $\phi_1, \phi_2, \ldots, \phi_K$. The selection of $K$ depends on the explained variance threshold.

For each initial condition and end time solution pair $u_i(x,\tau_{i-1}), u_i(x,\tau_i)$, we can represent them as weighted sum of principal components, say,
$u_i(x,\tau_{i-1}) = \sum_{j=1}^{K} \alpha_j(\tau_{i-1}) \phi_j(x), $ 
$ u_i(x,\tau_i) = \sum_{j=1}^{K} \beta_j(\tau_i) \phi_j(x).$
Our objective now becomes predicting the coefficients $\beta_j(\tau_i)$ of the end time solution using the coefficients $\alpha_j(\tau_{i-1})$ of the initial condition.

For illustration, we assume the FKE Eq.~(8) can be expressed by:
\begin{align}
\frac{\partial u_i}{\partial t}(x,t) 
&= \left( \mathcal{L} - \frac{1}{2} h^\top S^{-1} h \right) u_i(x,t)  \\
&= {} \; A(x, a(t)) \frac{\partial^2 u_i}{\partial x^2} 
+ B(x, b(t)) \frac{\partial u_i}{\partial x} \notag \\
&\quad + C(x, c(t)) u_i(x,t)  + D(x, d(t)) \notag \label{eq:time_variant_fke}
\end{align}
where $A$, $B$, $C$, and $D$ are coefficients of the terms on the right-hand side of the FKE, and $a(t)$, $b(t)$, $c(t)$, $d(t)$ are the time-dependent part of the coefficients.

Stage IB trains an approximate solver for time-variant FKE in a data-driven approach. The network inputs comprises two parts: the PCA-derived initial condition coefficients $\alpha_j(\tau_{i-1})$, $j = 1, 2, \ldots, K$, together with the time-dependent coefficients of the FKE $a(t), b(t), c(t), d(t)$ evaluated at $M+1$ time points $t = \tau_{i-1} + \frac{m}{M} \Delta t$, $m = 0, 1, 2, \ldots, M$. The network outputs predictions $\hat{\beta}_j(\tau_i)$ of the actual terminal coefficients $\beta_j(\tau_i)$, $j = 1, 2, \ldots, K$.
We note that for the time-invariant case, the only inputs are the $\alpha_j(\tau_{i-1})$, $j = 1, 2, \ldots, K$. The loss function is defined as the mean squared error between the predicted and actual coefficients:
\begin{equation}
L = \frac{1}{K} \sum_{j=1}^{K} \left[ \hat{\beta}_j(\tau_i) - \beta_j(\tau_i) \right]^2
\end{equation}

\begin{figure}[htbp]
    \centering
    \includegraphics[width=\linewidth]{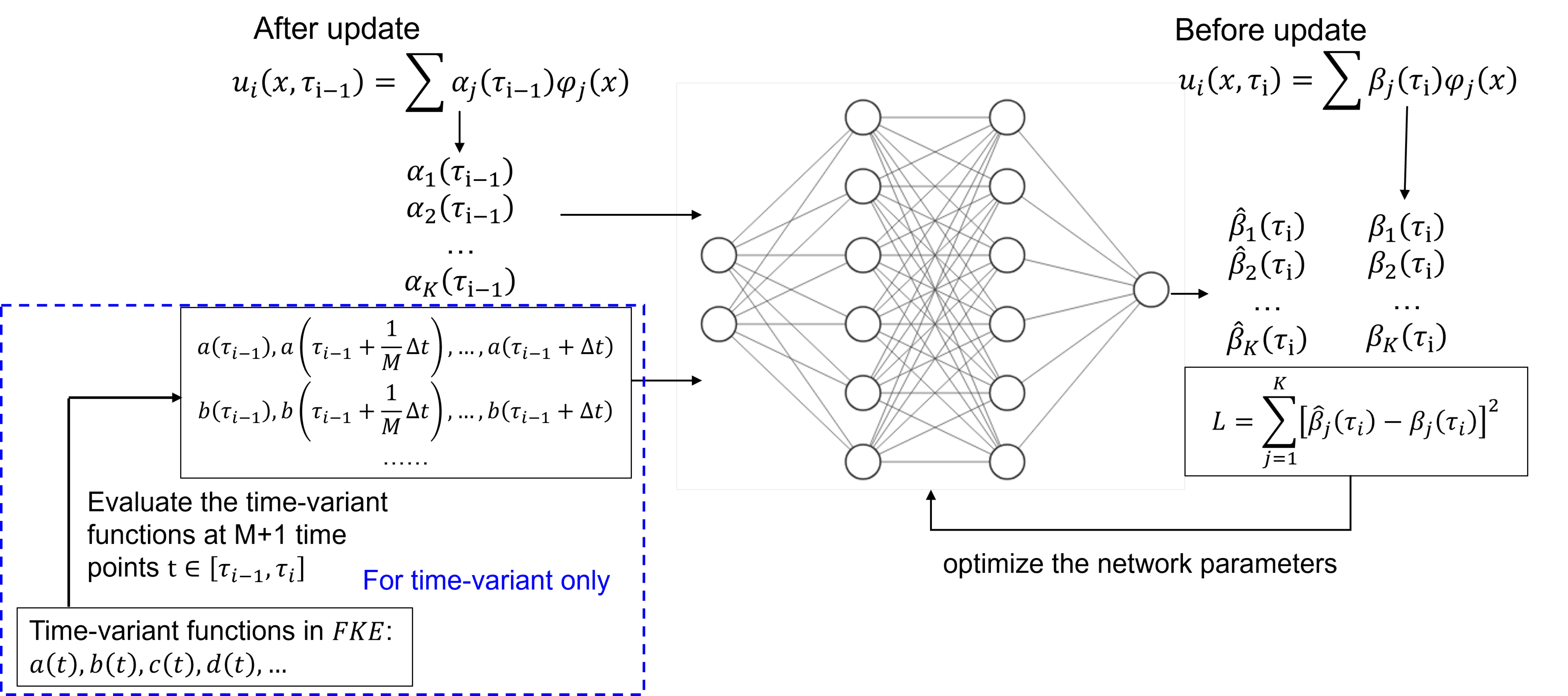}  
    \caption{The training of the approximate FKE solver based on PCA.}
    \label{fig:stageIB}
\end{figure}

A schematic of training process of this approximate FKE solver is depicted in Fig.~\ref{fig:stageIB}. Once trained, this approximate FKE solver can then be used in Stage II for fast prediction of end time solutions with very limited storage needs, enabling near-real-time online computation.

It is critical to ensure large and representative training dataset that includes sufficient pairs of initial conditions and solutions across various intervals. Such comprehensive training ensures: (1) the principal components can capture enough variability of $u$; (2) robustness to variations in time-variant coefficients; and (3) the accuracy of the approximate FKE solver.

\begin{algorithm}[htbp]
\caption{Train the approximate FKE solver (offline)}
\label{alg:stageIB}
\begin{algorithmic}[1]
\STATE \textbf{Initialization:} given $u_i(x, \tau_{i-1})$, $u_i(x, \tau_i)$, $K$, the time-variant coefficients in the FKE $a(t)$, $b(t)$, $\ldots$, the number of points $M$ to evaluate the coefficients, the maximum number of epochs $E_{\max}$, and the network structure, we initialize the neural network parameter $\eta$.
\STATE \textbf{PCA:} perform PCA on the $u_i(x, \tau_{i-1})$, $u_i(x, \tau_i)$, $i = 1,2,\ldots,N_T$, and extract the top $K$ principal components $\phi_1, \phi_2, \ldots, \phi_K$.
\STATE \textbf{Projection onto PCA Basis:} Decompose the $u_i(x, \tau_{i-1})$, $u_i(x, \tau_i)$ into weighted sum of principal components:
\[
u_i(x, \tau_{i-1}) = \sum_{j=1}^K \alpha_j(\tau_{i-1}) \phi_j(x), \quad
u_i(x, \tau_i) = \sum_{j=1}^K \beta_j(\tau_i) \phi_j(x).
\]
\FOR{$i = 1$ to $E_{\max}$}
    \STATE Feed the coefficients $\alpha_j(\tau_{i-1})$, $j = 1,2,\ldots,K$ and the time-variant functions $a(t)$, $b(t)$, $c(t)$, $d(t)$ evaluated at $M+1$ time points $t = \tau_{i-1} + \frac{m}{M} \Delta t$, $m = 0,1,2,\ldots,M$ into the network, and get output $\hat{\beta}_j(\tau_i)$, $j = 1,2,\ldots,K$.
    \STATE Compute the loss based on Eq.~(12).
    \STATE Optimize the network parameter $\eta$.
\ENDFOR
\STATE Store the trained model for real-time execution in Stage II.
\end{algorithmic}
\end{algorithm}

Both Algorithms \ref{alg:stageIA} and \ref{alg:stageIB} are computationally expensive. Their training should be done offline. The trained approximate FKE solver for time-variant system can then be applied to online execution in the following stage.

\subsection{Stage II (online): Real-time execution for system state estimation}

With the trained approximate FKE solver from Stage IB, we proceed to the online execution stage: the real-time computation of the system state with significantly reduced computation time and storage requirements. The essential components required for the online execution include: 

\begin{enumerate}
    \item the principal components given the principal components $\phi_j(x)$ from Stage IB, 
    \item the trained approximate FKE solver based on PCA in Stage IB, and 
    \item (for time-variant case only) the time-variant coefficients $a(t), b(t), \ldots$, and the number of points $M$ in each time interval for evaluation. 
\end{enumerate}

The primary storage requirement of this algorithm consists of the principal components, while the storage of network parameters requires little storage. Compared to the approach taken by previous method \cite{Dong2020}, the overall storage is relatively small because only one universal approximate FKE solver is stored rather than multiple solvers for different time intervals. 

We estimate the system state recursively. For each time step $[\tau_{i-1}, \tau_i]$, we execute the algorithm as follows:

\begin{algorithm}[htbp]
\caption{Real-time computation of the system state (online)}
\label{alg:stage_ii}
\begin{algorithmic}[1]
\STATE Given the principal components $\phi_j(x)$, the time-variant coefficients $a(t), b(t), \ldots$, and the number of points $M$ in each time interval for evaluation, the approximate FKE solver based on PCA in Stage IB, the initial distribution function $\sigma_0(x)$, number of time steps $N_T$, and the observation signals $y_{\tau_i}$ for $i = 0,1,2,\ldots,N_T$. We initialize the initial probability distribution $u_0(x,0) = \sigma_0(x)$.
\FOR{$i = 1$ to $N_T$}
    \STATE Accept signals $y_{\tau_{i-1}}$ and update and obtain initial solution $u_i(x, \tau_{i-1})$
    \STATE Obtain the PC coefficients $\alpha_j(\tau_{i-1})$ of the initial condition
    \STATE Obtain end time PC coefficients $\beta_j(\tau_i)$ using the time-variant approximate FKE solver with the time-variant coefficients
    \STATE Obtain the end time solution 
    \[
    u_i(x, \tau_i) = \sum_{j=1}^{K} \beta_j(\tau_i) \phi_j(x)
    \]
\ENDFOR
\end{algorithmic}
\end{algorithm}

\begin{figure}[htbp]
    \centering
    \includegraphics[width=\linewidth]{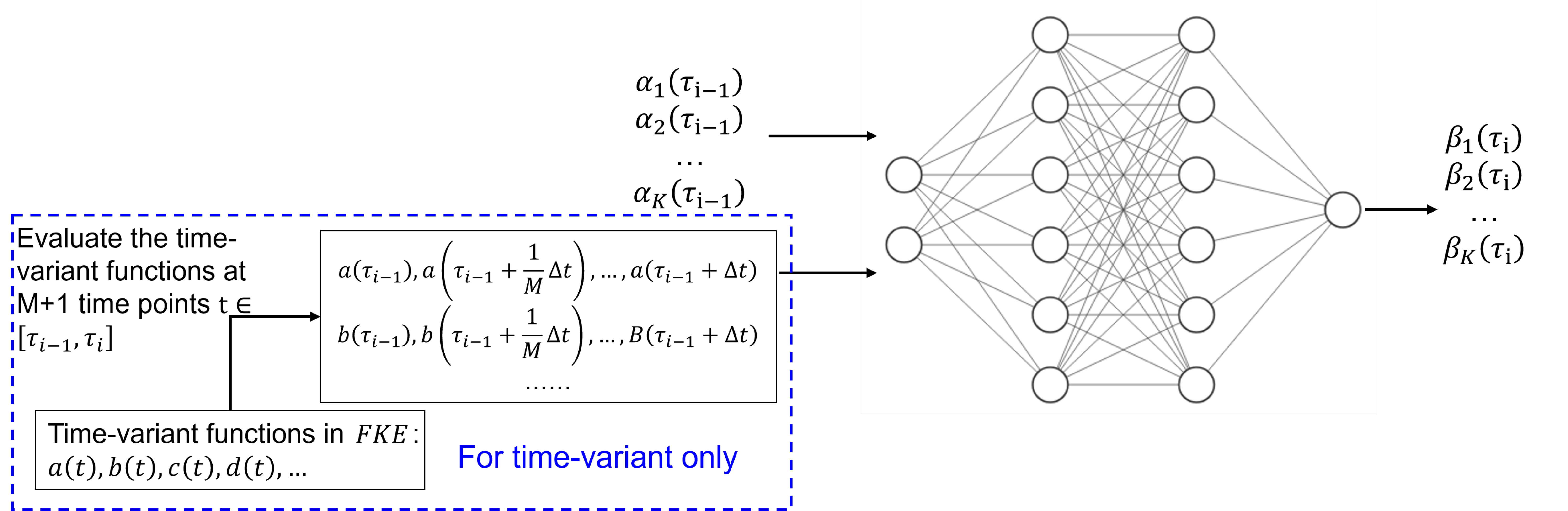}  
    \caption{Execute the approximate FKE solver and obtain the end time solutions.}
    \label{fig:stage2}
\end{figure}

The combination of Algorithms 1, 2, and 3 is referred to as the PINN-Based Yau-Yau Filter, abbreviated as PINNYYF.

\begin{remark}
 Our proposed algorithm constitutes the first global Yau-Yau filter incorporating data-driven adaptive basis functions. It retains the global characteristics of the original Yau-Yau method while also leveraging the efficiency of data-driven methods. Compared to RNN-YYF~\cite{Chen2024} and neural projection filters~\cite{Tao2023}, our algorithm is capable of reconstructing the probability density function. Furthermore, compared to existing PINN-based Yau-Yau filters~\cite{Shi2024}, our method provides adaptive and optimal set of basis functions, overcoming the limitations of traditional Galerkin methods in basis function selection.    
\end{remark}

\begin{remark}
Assuming the network consists of \( L \) layers, with \( h \) nodes in each layer, the total number of network parameters \( P \) is approximately \( P \approx Lh^2 \). Assuming the PDE is of order \( k \), and there are \( N_{\text{PDE}}, N_{\text{IC}}, N_{\text{BC}} \) collocation points for the PDE, initial condition, and boundary condition, respectively, the total feed-forward complexity for the training in one epoch is 
\[
\mathcal{O}((kN_{\text{PDE}} + N_{\text{IC}} + N_{\text{BC}}) L h^2).
\]
Assuming the network is trained for \( E \) epochs, the total complexity for one time step is 
\[
\mathcal{O}((kN_{\text{PDE}} + N_{\text{IC}} + N_{\text{BC}}) E L h^2).
\]
When the hyperparameters are adjusted properly, each time step will only require a few epochs to train, resulting in decent time efficiency. 
\end{remark}

\section{Numerical Experiments}

\subsection{Network Description}

We apply the proposed PINNYYF to three examples and test their efficiencies in terms of accuracy, time cost, and storage cost. The first example focuses on a time-invariant system, and the second and third examples involve time-variant systems. We also implement the EKF and PF as baseline for comparison across all examples to evaluate our method’s performance against traditional nonlinear filter techniques.

In Stage IA, the PINN-based FKE solver is configured as a fully connected neural network. Key hyperparameters for this stage are summarized in Table~\ref{tab:stage1a_hyperparams}. The convergence threshold $\epsilon$ is adjusted for each specific example to ensure the solution meets the required accuracy.

\begin{table}[htbp]
\centering
\caption{Hyperparameter settings for Stage IA.}
\label{tab:stage1a_hyperparams}
\begin{tabular}{l c}
\toprule
\textbf{Parameters} & \textbf{Values} \\
\midrule
Hidden layers                  & 4 \\
Nodes per layer               & 40 \\
Activation function           & Tanh \\
Optimizer                     & Adam \cite{Kingma2014} \\
Learning rate                 & 0.001 \\
Maximum epochs each time step & 10000 \\
\bottomrule
\end{tabular}
\end{table}

In Stage IB, recognizing the solution $u_i(x, \tau_i)$ at each time step differs only slightly from the initial condition $u_i(x, \tau_{i-1})$, particularly for small $\Delta t$, we adopt a residual learning approach. Instead of directly learning the terminal coefficients $\beta_j(\tau_i)$, we predict the difference between the initial principal components (PC) coefficients $\alpha_j(\tau_{i-1})$ and the end time coefficients $\beta_j(\tau_i)$. This approach simplifies the learning process and improves training efficiency.

\begin{figure}[htbp]
    \centering
    \includegraphics[width=\linewidth]{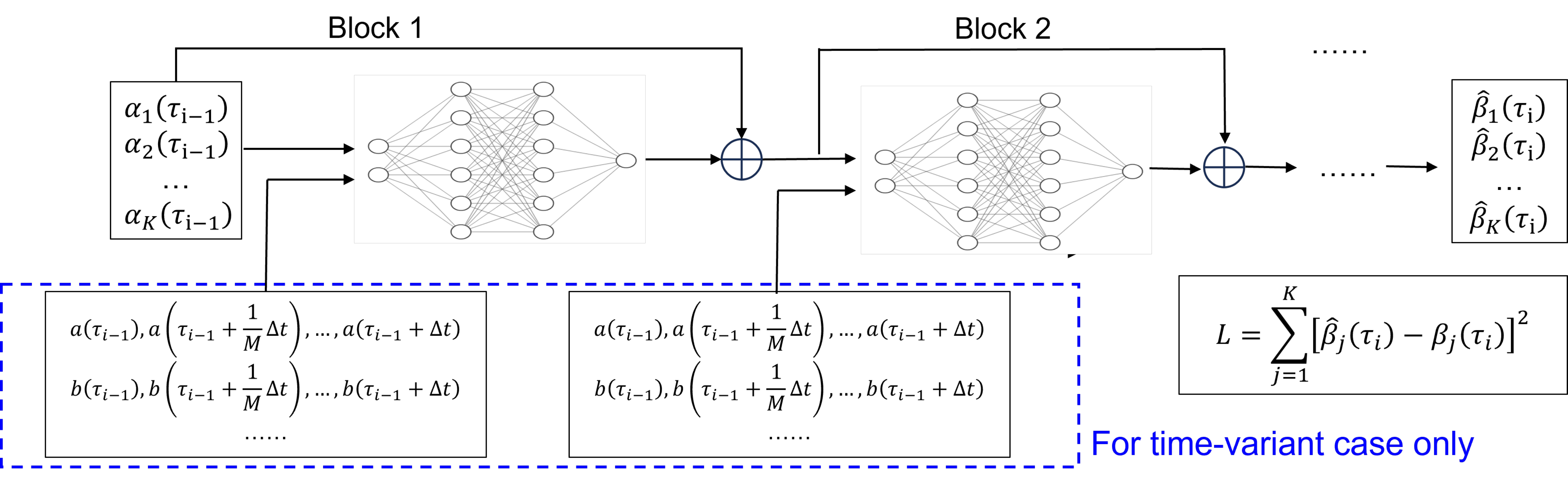} 
    \caption{The residual network structure used in Stage IB.}
    \label{fig:resnet_structure}
\end{figure}

To achieve this, we adopt the residual network structure~\cite{He2015} as the approximate FKE solver network. The residual network consists of multiple residual blocks, each of which consists of a fully connected network with a skip connection, as depicted in Fig.~\ref{fig:resnet_structure}. The skip connection adds the network’s input directly to its output for each block, improving the gradient flow and enhancing training stability.

The input to each residual block comes from two parts:
\begin{enumerate}
    \item The output coming from the preceding block (or the initial PC coefficient $\alpha_j(\tau_{i-1})$ for the first block).
    \item (Only for time-variant cases) Time-dependent coefficients $a(t), b(t), c(t), d(t), \ldots$ evaluated at $M+1$ time points $t = \tau_{i-1} + \frac{m}{M} \Delta t$, $m = 0, 1, 2, \ldots, M$ within each interval.
\end{enumerate}

For consistency across all experiments, we configure the residual network with 3 residual blocks, each is a fully connected network with 2 layers and 64 nodes per layer. Training uses the Adam~\cite{Kingma2014} optimizer with a learning rate of 0.001 over 20,000 epochs for convergence.

We measure the accuracy of each algorithm using mean squared error (MSE) for each state variable, defined as:
\begin{equation}
\text{MSE}_{x_l} = \frac{1}{N_T + 1} \sum_{i=0}^{N_T} \left( x_l(\tau_i) - \hat{x}_l(\tau_i) \right)^2
\end{equation}
where $x_l(\tau_i)$ is the true value of the $l$-th component of the state variable at time $\tau_i$ and $\hat{x}_l(\tau_i)$ is the corresponding estimated value, for $l = 1, 2, \ldots, d$.

All numerical experiments were conducted on an NVIDIA RTX3070 GPU and a platform with 16 Intel Core i7-10700 CPUs at 2.90~GHz. The neural networks were implemented using PyTorch.

\subsection{Example 1: Cubic Sensor Problem (Time-invariant)}

In the first example, we demonstrate the effectiveness of PINNYYF in strongly nonlinear time-invariant system using the two-dimensional cubic sensor problem:
\begin{equation}
\begin{cases}
dx_1(t) = [-0.4 x_1(t) + 0.1 x_2(t)]\,dt + dv_1(t) \\
dx_2(t) = -0.6 x_2(t)\,dt + dv_2(t) \\
dy_1(t) = x_1^3(t)\,dt + dw_1(t) \\
dy_2(t) = x_2^3(t)\,dt + dw_2(t)
\end{cases}
\end{equation}

where $v(t) = [v_1(t), v_2(t)]^\top$, $w(t) = [w_1(t), w_2(t)]^\top$ are Brownian motion processes with 
\[
\mathbb{E}(dv(t) dv(t)^\top) = I_2\,dt,\quad
\mathbb{E}(dw(t) dw(t)^\top) = I_2\,dt,
\]
where $I_2$ is the 2-dimensional identity matrix. The time interval is $\Delta t = 0.01$, the total number of time steps is $N_T = 5000$, and the initial state $[x_1(0), x_2(0)]$ is sampled using the normal distribution $\mathcal{N}((0, 0), 0.2 I_2)$.

The FKE equation corresponding to Eq.~(14) is:
\begin{align}
\frac{\partial u}{\partial t} 
&= \left( \mathcal{L} - \frac{1}{2} h^\top S^{-1} h \right) u \notag \\
&= \frac{1}{2} \frac{\partial^2 u}{\partial x_1^2}
+ \frac{1}{2} \frac{\partial^2 u}{\partial x_2^2}
+ u \notag\\
&\quad - [-0.4 x_1 + 0.1 x_2] \frac{\partial u}{\partial x_1}
+ 0.6 x_2 \frac{\partial u}{\partial x_2} \notag\\
 & \quad - \frac{1}{2} (x_1^6 + x_2^6) u 
\label{eq:fke_cubic_sensor}
\end{align}

The relationship for absorbing observations and updating initial conditions becomes:
\begin{equation}
u_i(x, \tau_i) = \exp\left\{ 
(x_1^3,x_2^3)^\top 
\left( y_{\tau_i} - y_{\tau_{i-1}} \right)
\right\}
\cdot u_{i-1}(x, \tau_i). 
\label{eq:update_obs}
\end{equation}

We assume the initial distribution for the first time step is:
$u_0(x, 0) = \exp\left( -\frac{1}{2} |x|^2 \right)$.

In Stage IA, we train the PINN-based FKE solver with computational domain $[-2.2, 2.2] \times [-2.2, 2.2]$ and convergence threshold $\epsilon = 0.001$, and conduct 20 Monte Carlo simulations. The average MSE values for the 20 simulations are $\text{MSE}_{x_1} = 0.419$ and $\text{MSE}_{x_2} = 0.407$. The mean CPU time for the 20 simulations was 3.5 hours, with an average of 298 epochs per time step.

\begin{figure}[htbp]
    \centering
    \begin{subfigure}[t]{0.48\linewidth}
        \centering
        \includegraphics[width=\linewidth]{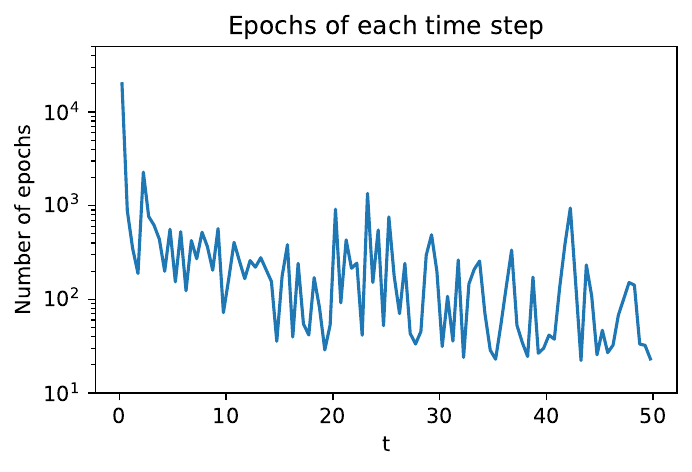}  
        \caption{}
    \end{subfigure}
    \hfill
    \begin{subfigure}[t]{0.48\linewidth}
        \centering
        \includegraphics[width=\linewidth]{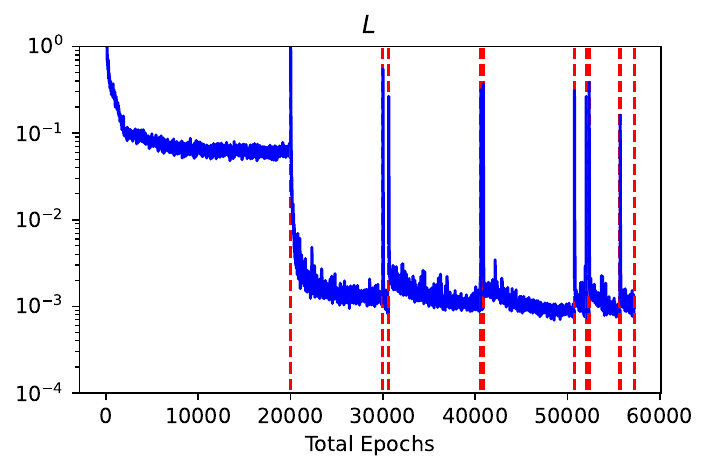}  
        \caption{}
    \end{subfigure}
    \caption{(a) The number of epochs for each time step; (b) the loss curves for the time step 1 to 10.Each time step is
represented by the interval between two vertical lines.}
    \label{fig:epochs_and_loss}
\end{figure}

Fig.~\ref{fig:epochs_and_loss}(a) illustrates the number of epochs required for each time step, revealing that the number of epochs decreases significantly after a few initial time steps. This suggests that the network adapts to different initial conditions quickly after a few time steps. Fig.~\ref{fig:epochs_and_loss}(b) depicts the loss curves for the first 10 time steps, where the number of epochs is the highest at the first step and decreases quickly in the subsequent steps, and the loss drops quickly after a few epochs in subsequent steps. This rapid convergence after a few initial training steps suggests that the PINN is efficient in capturing the system's dynamics, leading to reduced training epochs for later time steps.

\begin{figure}[htbp]
    \centering
    \includegraphics[width=\linewidth]{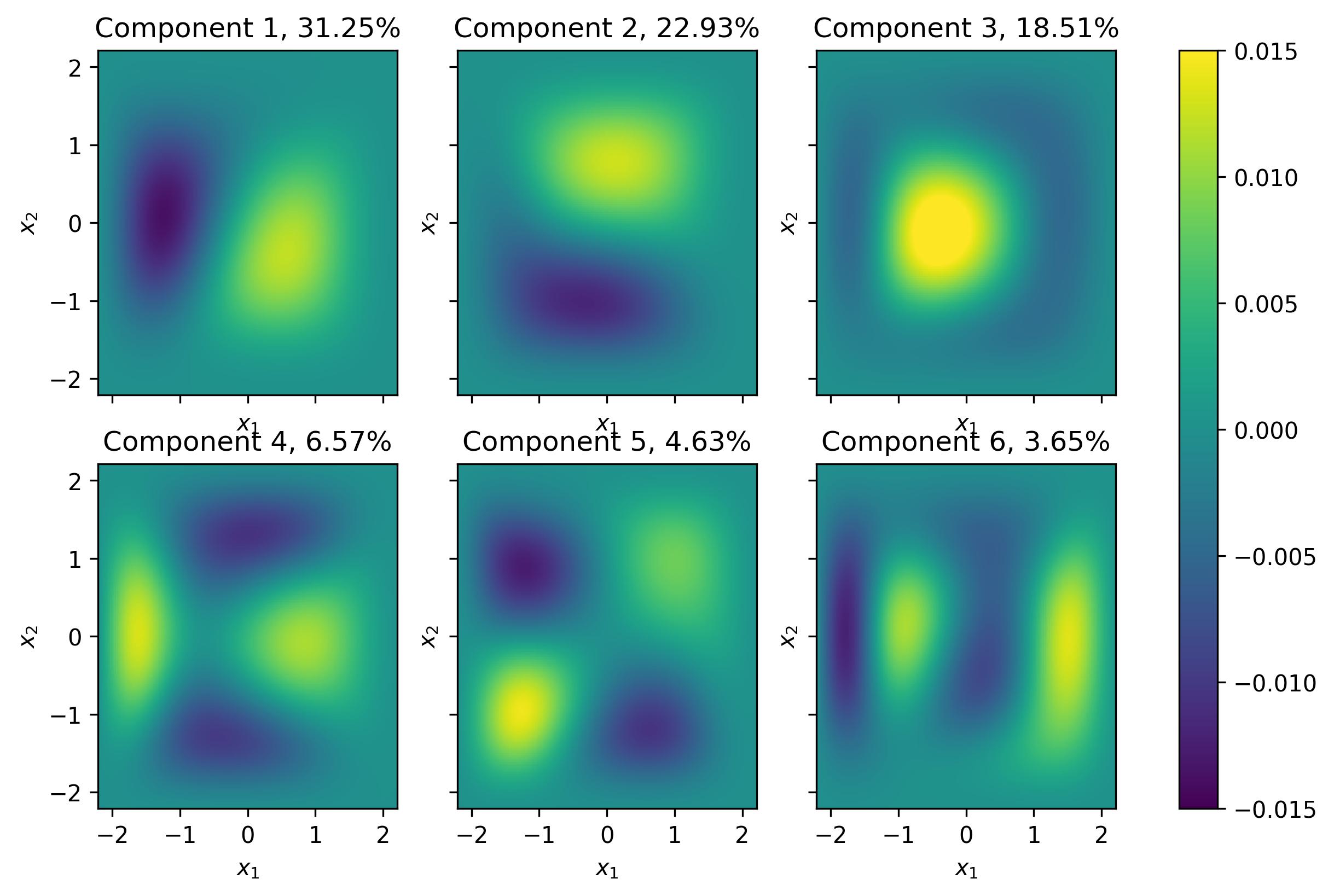}  
    \caption{First 6 principal components and their explained variance ratios.}
    \label{fig:pca_components}
\end{figure}

In Stage IB, we extract 30 principal components that account for 99.7\% of the explained variance. Fig.~\ref{fig:pca_components} visualizes the first 6 principal components and their explained variance ratios. The high variance captured suggests that these components effectively summarize the system’s dynamics. We train the approximate FKE solver network with 32,000 time steps as training set and 8,000 time steps as testing set with a batch size $B = 512$ for training.

In Stage II, we employed the network trained in Stage IB to predict real-time state variables $x_1, x_2$ based on real-time observation signals $y_{\tau_i}$. We compare PINNYYF with EKF, PF, and the Yau-Yau filter based on Legendre spectral method (LSMYYF)~\cite{Dong2020}. The PF is employed with 100 particles. Table~\ref{tab:eg1MSE} summarizes the comparison between the three methods in terms of accuracy (MSE), CPU time, and storage efficiency.

The EKF method fails in the accuracy mainly due to the high nonlinearity. The PF method achieves decent accuracy but is very time consuming. The spectral method achieves satisfactory results with decent accuracy, showing the spectral method is good at dealing with time-invariant problems. Our method strikes a balance, achieving comparable accuracy with LSMYYF while offering moderate CPU time and storage efficiency. Moreover, the MSE for PINNYYF is only slightly larger than those from PINN (Stage IA offline results), but with much shorter time cost. This shows our dimension reduction and approximate solver for FKE is effective and efficient.

\begin{figure}[htbp]
    \centering
    \includegraphics[width=\linewidth]{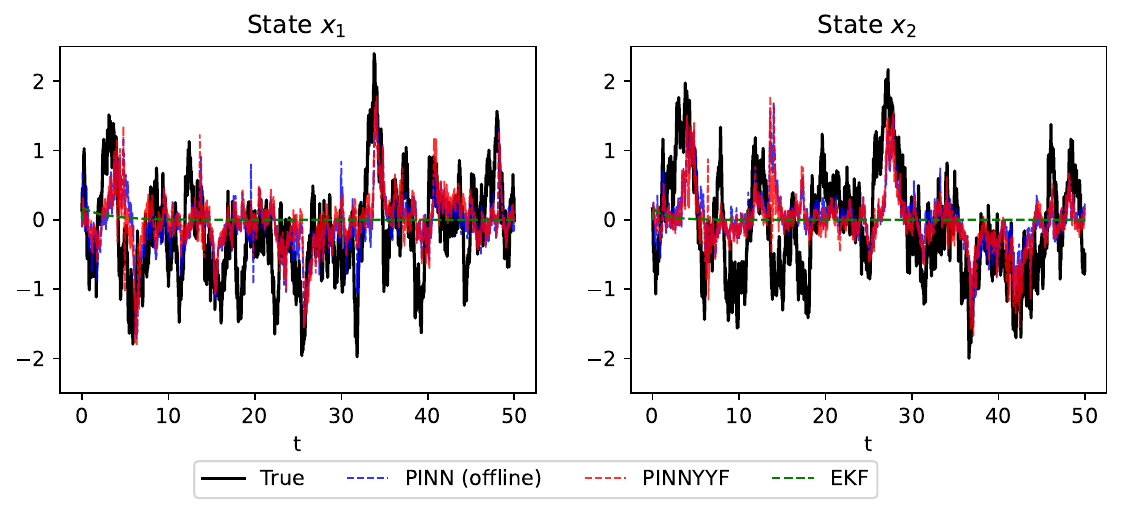}
    \caption{The recovery of state variables $x_1$, $x_2$ using PINN (offline), PINNYYF, and EKF, compared to true values.}
    \label{fig:eg1state}
\end{figure}

Fig.~\ref{fig:eg1state} shows a representative realization of the state variables $x_1, x_2$, from the 20 Monte Carlo simulations. The curves show the true state variables alongside the predictions from the EKF, PF, PINN (offline), and PINNYYF. Our method tracks the states closely over time, whereas the EKF fails to capture the full variation, particularly in regions with high nonlinearity. This demonstrates the robustness of PINNYYF in accurately tracking nonlinear state dynamics in real time.

\begin{figure}[htbp]
    \centering
    \includegraphics[width=\linewidth]{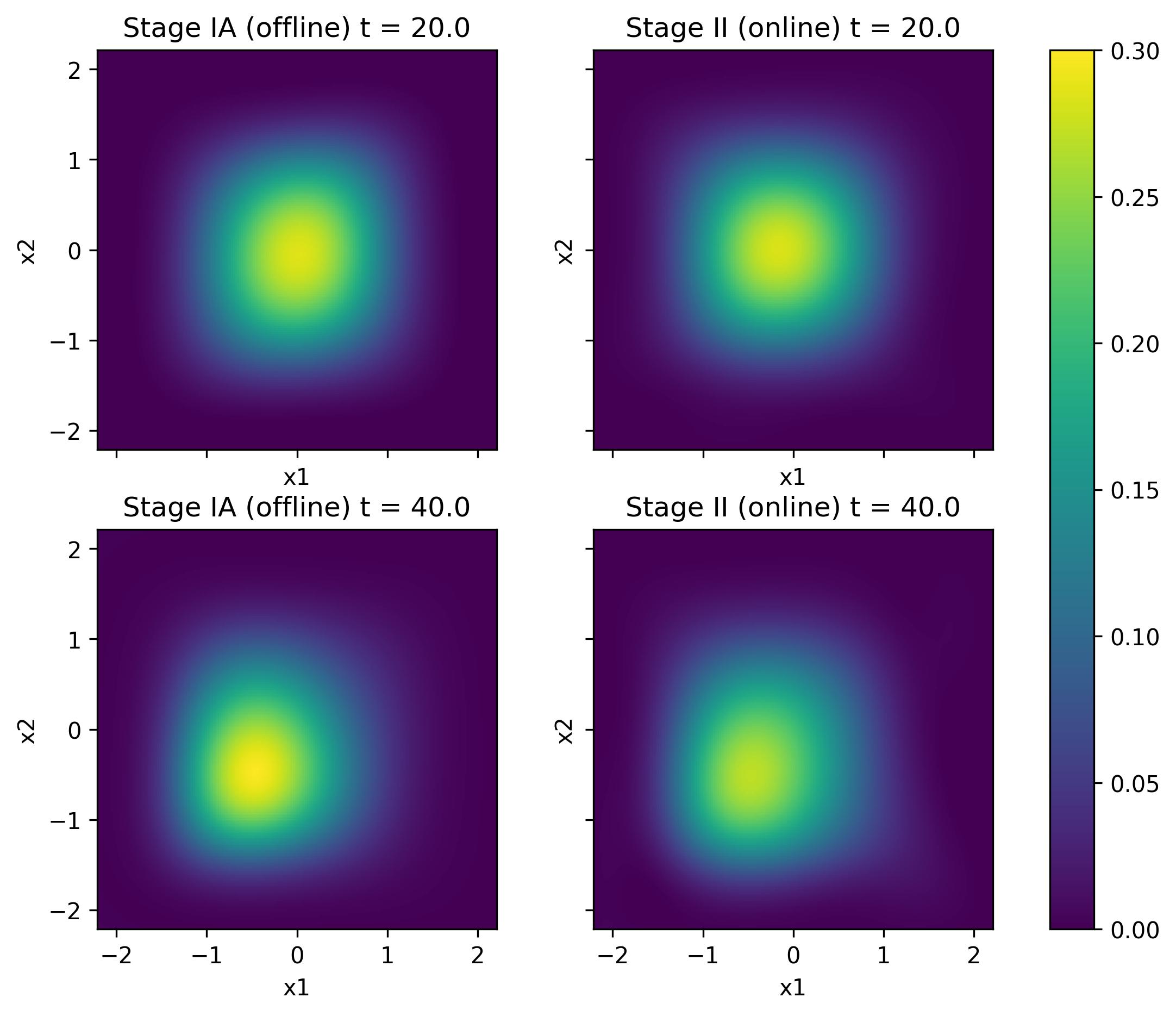}  
    \caption{Two different snapshots for the probability density obatined $u$ in stage IA (left column) and the reconstructed probability density in stage II (right column).}
    \label{fig:eg1pdf}
\end{figure}

Fig.~\ref{fig:eg1pdf} illustrates the probability density obtained during Stage IA, alongside the reconstructed probability density from Stage II, at two selected time snapshots. This comparison demonstrates that the reconstructed probability density in the online stage closely aligns with the one obtained in the offline stage, thereby validating the effectiveness of our proposed algorithm.

\begin{table}[htbp]
\centering
\caption{Example 1: Comparison of accuracy, CPU time per time step, and storage for EKF, LSMYYF, PINN (offline), and PINNYYF.}
\label{tab:eg1MSE}
\resizebox{\linewidth}{!}{
\begin{tabular}{lcccc}
\toprule
\textbf{Methods} & \textbf{MSE-$x_1$} & \textbf{MSE-$x_2$} & \textbf{CPU time} & \textbf{Storage} \\
\midrule
EKF                & 1.013 & 0.759 & 0.06\,ms   & $<$1\,kb \\
PF (100 particles) & 0.435 & 0.423 & 3.68\,ms   & $<$1\,kb \\
LSMYYF ($20^2$ basis) & 0.409 & 0.424 & 0.77\,ms   & 1,251\,kb \\
PINN (offline)     & 0.419 & 0.407 & 2.52\,s    & 318\,kb \\
PINNYYF            & 0.432 & 0.420 & 0.78\,ms   & 1,159\,kb \\
\bottomrule
\end{tabular}%
}
\end{table}

\subsection{Example 2: Time-variant Harmonic Sensor Problem}

Example 2 is a time-variant weakly nonlinear problem selected from the numerical examples from Dong et al.~\cite{Dong2020}. The system dynamics are given by:
\begin{equation}
\begin{cases}
dx_1(t) = [-0.4 x_1(t) + 0.1 x_2(t)]\,dt + [1 + 0.1 \cos(20\pi t)]\,dv_1(t) \\
dx_2(t) = -0.6 x_2(t)\,dt + [0.9 + 0.2 \cos(18\pi t)]\,dv_2(t) \\
dy_1(t) = x_1(t)\,[1 + 0.2 \cos(x_2(t))]\,dt + dw_1(t) \\
dy_2(t) = x_2(t)\,[1 + 0.2 \cos(x_1(t))]\,dt + dw_2(t)
\end{cases}
\label{eq:example2_system}
\end{equation}

where $v(t) = [v_1(t), v_2(t)]^\top$, $w(t) = [w_1(t), w_2(t)]^\top$ are Brownian motion processes with
\[
\mathbb{E}[dv(t) dv(t)^\top] = I_2\,dt, \quad \mathbb{E}[dw(t) dw(t)^\top] = I_2\,dt,
\]
and $I_2$ is the 2-dimensional identity matrix. 

Here, we let the time interval be $\Delta t = 0.01$, the number of time steps be $N_T = 5000$, and the initial state $[x_1(0), x_2(0)]$ is sampled from the normal distribution $\mathcal{N}((0,0), 0.2 I_2)$.

The FKE equation corresponding to Eq.~(16) is:
\begin{align}
\frac{\partial u}{\partial t} &= \left( \mathcal{L} - \frac{1}{2} h^\top S^{-1} h \right) u \notag \\
&= \frac{1}{2} \left[1 + 0.1 \cos(20\pi t)\right]^2 \frac{\partial^2 u}{\partial x_1^2}
+ \frac{1}{2} \left[0.9 + 0.2 \cos(18\pi t)\right]^2 \frac{\partial^2 u}{\partial x_2^2} \notag \\
&\quad + u
- [-0.4 x_1 + 0.1 x_2] \frac{\partial u}{\partial x_1}
+ 0.6 x_2 \frac{\partial u}{\partial x_2} \notag \\
&\quad - \frac{1}{2} x_1^2 \left[1 + 0.2 \cos(x_2)\right]^2 u
- \frac{1}{2} x_2^2 \left[1 + 0.2 \cos(x_1)\right]^2 u
\label{eq:fke_harmonic}
\end{align}

The updated initial condition for the $i$-th time interval at $t = \tau_i$ is:
\begin{equation}
u_i(x, \tau_i) = \exp \left\{
\begin{pmatrix}
x_1 [1 + 0.2 \cos(x_2)] \\
x_2 [1 + 0.2 \cos(x_1)]
\end{pmatrix}^\top
(y_{\tau_i} - y_{\tau_{i-1}})
\right\}
\cdot u_{i-1}(x, \tau_i)
\label{eq:obs_update_harmonic}
\end{equation}

The initial distribution for $t = \tau_0 = 0$ is $u_0(x, 0) = \exp\left( -\frac{1}{2} |x|^2 \right)$.

In Stage IA, the computational domain for the solution of FKE is $[-3, 3] \times [-3, 3]$, and the convergence threshold is $\epsilon = 0.0002$ for the FKE solver training. The MSEs for 20 Monte Carlo simulations are $\text{MSE}_{x_1} = 0.602$ and $\text{MSE}_{x_2} = 0.478$. The mean CPU time for the 20 simulations is 1.9 hours, and the average number of epochs for each time step is 85.

In Stage IB, we extract 30 principal components which account for 99.97\% of the explained variance. The approximate FKE solver inputs the first 30 PC coefficients, and the time-variant functions $\cos(20\pi t)$ and $\cos(18\pi t)$ evaluated at $M = 4$ time points between $\tau_{i-1}$ and $\tau_i$. We train the approximate FKE solver network with 32,000 time steps as training set, and 8,000 time steps as testing set, with a batch size $B = 512$ for the training.

We compare our implementation method with EKF, PF (100 particles), and LSMYYF. Table~\ref{tab:eg2MSE} lists the MSEs of state variables for the three filter methods. Due to the weak nonlinearity, the EKF performs reasonably well. PINNYYF exhibits comparable accuracy to EKF, outperforming PF in both accuracy and CPU time. Compared to spectral-methods-based Yau-Yau filter, PINNYYF also achieves better accuracy with much less storage cost than LSMYYF.

\begin{figure}[htbp]
    \centering
    \includegraphics[width=\linewidth]{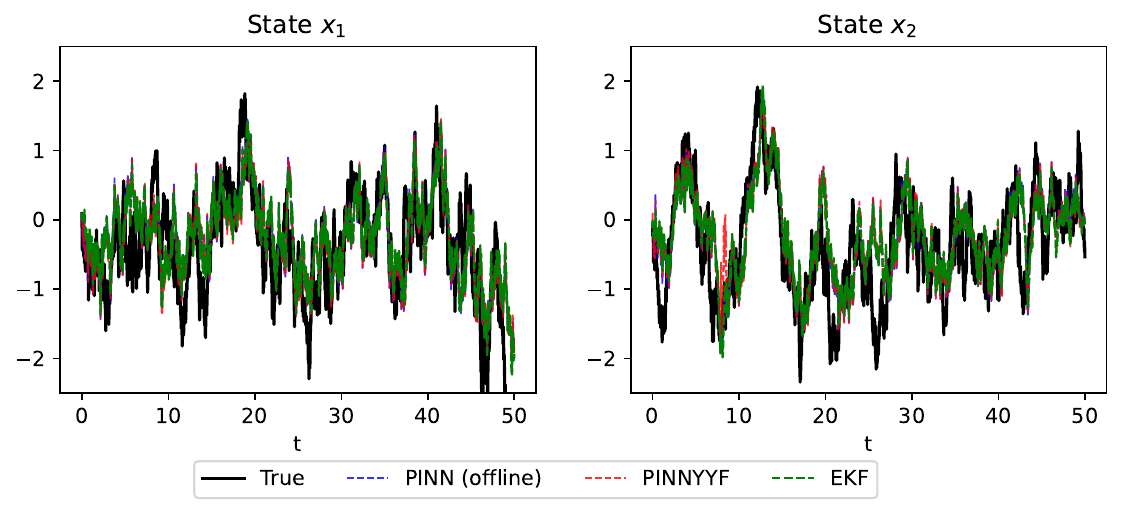}  
    \caption{The recovery of state variables $x_1$ and $x_2$ using PINN (offline), PINNYYF, and EKF, compared with true values.}
    \label{fig:eg2state}
\end{figure}

Fig.~\ref{fig:eg2state} shows one of the 20 simulations of the state variables $x_1$ and $x_2$. The figure compares the variable curves obtained using these three methods with the true curve, including the EKF, spectral method and PINYYF. We see that PINNYYF performs well tracking the states throughout time.

\begin{table}[htbp]
\centering
\caption{Example 2: Comparison of accuracy, CPU time per time step and storage for EKF, PF, LSMYYF, PINN (offline), and PINNYYF.}
\label{tab:eg2MSE}
\resizebox{\linewidth}{!}{%
\begin{tabular}{lcccc}
\toprule
\textbf{Methods} & \textbf{MSE-$x_1$} & \textbf{MSE-$x_2$} & \textbf{CPU time} & \textbf{Storage} \\
\midrule
EKF                   & 0.611 & 0.468 & 0.06\,ms   & $<$1\,kb \\
PF (100 particles)    & 0.631 & 0.489 & 3.86\,ms   & $<$1\,kb \\
LSMYYF ($10^2$ basis) & 0.634 & 0.498 & 0.18\,ms   & 390,626\,kb \\
PINN (offline)        & 0.602 & 0.478 & 1.36\,s    & 318\,kb \\
PINNYYF               & 0.621 & 0.488 & 0.84\,ms   & 1,178\,kb \\
\bottomrule
\end{tabular}%
}
\end{table}

\subsection{Example 3: Time-variant Cubic Sensor Problem}

This example is an adaptation of Example 1, with time-dependence in the state equations and high nonlinearity in the observation equations:

\begin{equation}
\begin{cases}
dx_1(t) = [-0.4 x_1(t) + 0.1 x_2(t)]\,dt + [1 + 0.1 \cos(20\pi t)]\,dv_1(t) \\
dx_2(t) = -0.6 x_2(t)\,dt + [0.9 + 0.2 \cos(18\pi t)]\,dv_2(t) \\
dy_1(t) = x_1^3(t)\,dt + dw_1(t) \\
dy_2(t) = x_2^3(t)\,dt + dw_2(t)
\end{cases}
\label{eq:example3_system}
\end{equation}

where $v(t) = [v_1(t), v_2(t)]^\top$, $w(t) = [w_1(t), w_2(t)]^\top$ are Brownian motion processes with
\[
\mathbb{E}[dv(t)\,dv(t)^\top] = I_2\,dt, \quad \mathbb{E}[dw(t)\,dw(t)^\top] = I_2\,dt,
\]
and $I_2$ is the 2-dimensional identity matrix. We set the time interval $\Delta t = 0.01$, total number of time steps $N_T = 5000$, and the initial state $[x_1(0), x_2(0)]$ is sampled from the normal distribution $\mathcal{N}((0,0), 0.2 I_2)$.

The FKE equation corresponding to Eq.~\eqref{eq:example3_system} is:
\begin{align}
\frac{\partial u}{\partial t} &= \left( \mathcal{L} - \frac{1}{2} h^\top S^{-1} h \right) u \notag \\
&= \frac{1}{2} \left[1 + 0.1 \cos(20\pi t)\right]^2 \frac{\partial^2 u}{\partial x_1^2}  \notag \\
&+ \frac{1}{2} \left[0.9 + 0.2 \cos(18\pi t)\right]^2 \frac{\partial^2 u}{\partial x_2^2} \notag \\
&+ u
- [-0.4 x_1 + 0.1 x_2] \frac{\partial u}{\partial x_1}
+ 0.6 x_2 \frac{\partial u}{\partial x_2} \notag\\
&- \frac{1}{2} \left( x_1^6 + x_2^6 \right) u
\label{eq:fke_example3}
\end{align}

The updated initial condition is:
\begin{equation}
u_i(x, \tau_i) = \exp\left\{
(x_1^3, x_2^3)^\top
(y_{\tau_i} - y_{\tau_{i-1}})
\right\}
\cdot u_{i-1}(x, \tau_i), \quad i \geq 2
\label{eq:update_example3}
\end{equation}

The initial distribution at $t = \tau_0 = 0$ is given by $u_0(x, 0) = \exp\left( -\frac{1}{2} |x|^2 \right)$.

In Stage IA, the computational domain for solving the FKE is $[-2.2, 2.2] \times [-2.2, 2.2]$, and the convergence threshold is $\epsilon = 0.001$. The MSEs for 20 Monte Carlo simulations are $\text{MSE}_{x_1} = 0.442$ and $\text{MSE}_{x_2} = 0.396$. The mean CPU time for the 20 simulations is 4.0 hours, and the average number of epochs for each time step is 302.

In Stage IB, 30 principal components are extracted which account for 99.6\% of the variance. The approximate FKE solver inputs 30 PC coefficients, and the time-variant functions $\cos(20\pi t)$ and $\cos(18\pi t)$ evaluated at 4 time points between $\tau_{i-1}$ and $\tau_i$, i.e., inputs of 38-dimensional vectors. We train the approximate FKE solver network with 32,000 time steps as training set, and 8,000 time steps as testing set with batch size $B = 512$.

Table~\ref{tab:eg3MSE} lists the MSEs of two state variables for our method and the EKF and PF methods. PINNYYF shows better accuracy with acceptable time-efficiency and storage-efficiency. Compared to the results of using PINN as FKE solver in Stage IA, the PINNYYF greatly reduces the CPU time with only a slight increase in MSE. PINNYYF also achieves better accuracy than PF with less CPU time.

\begin{figure}[htbp]
    \centering
    \includegraphics[width=\linewidth]{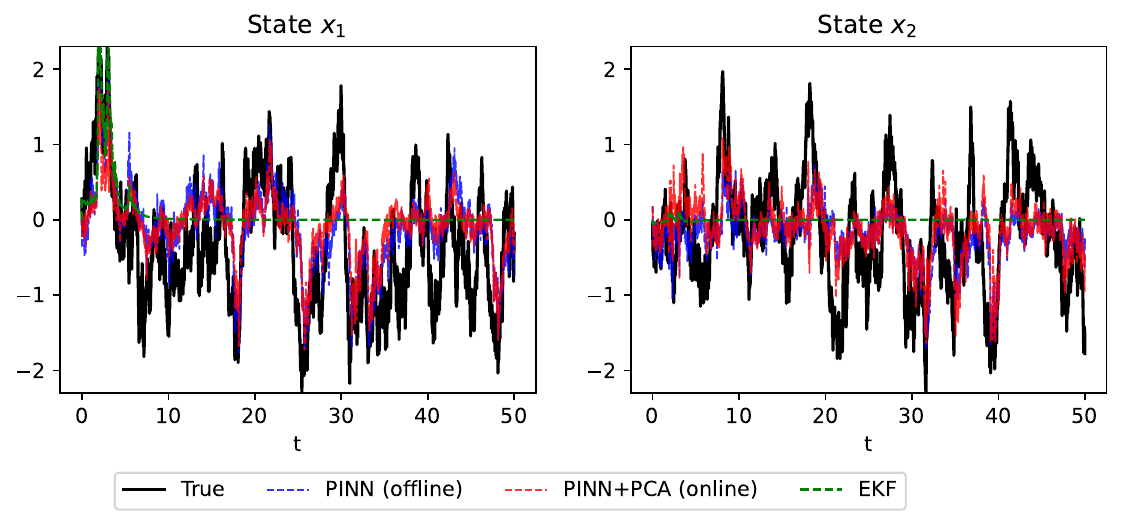}  
    \caption{The recovery of state variables $x_1$, $x_2$ using PINN (offline), PINNYYF, and EKF, compared to true values.}
    \label{fig:eg3state}
\end{figure}

Fig.~\ref{fig:eg3state} shows the state variables $x_1$ and $x_2$ from one of the 20 Monte Carlo simulations. We see that compared to EKF, PINNYYF closely captures the true state dynamics throughout time.

\begin{table}[htbp]
\centering
\caption{Example 3: Comparison of accuracy, CPU time per time step and storage for EKF, PINN (offline), and PINNYYF.}
\label{tab:eg3MSE}
\resizebox{\linewidth}{!}{%
\begin{tabular}{lcccc}
\toprule
\textbf{Methods} & \textbf{MSE-$x_1$} & \textbf{MSE-$x_2$} & \textbf{CPU time} & \textbf{Storage} \\
\midrule
EKF              & 0.913 & 0.759 & 0.06\,ms   & $<$1\,kb \\
PF (100 particles) & 0.455 & 0.415 & 3.87\,ms   & $<$1\,kb \\
PINN (offline)   & 0.442 & 0.396 & 2.88\,s    & 318\,kb \\
PINNYYF          & 0.453 & 0.405 & 0.86\,ms   & 1,319\,kb \\
\bottomrule
\end{tabular}%
}
\end{table}

\section{Conclusion}
This study presented a novel algorithm, PINNYYF, integrating PINNs, PCA, and Yau-Yau filtering framework, for addressing challenges posed by time-variant nonlinear filtering. The methodology utilizes PINNs to solve the FKE in the offline phase (Stage IA), and subsequently incorporates PCA to effectively compress solution data into PC coefficients, facilitating the creation of an approximate FKE solver (Stage IB). The solver is adapted to handle time-dependent problems by incorporating time-variant coefficients into the network's inputs. In the subsequent online phase (Stage II), this approximate solver enables near-real-time computation, requiring only the determination of PCA coefficients and updating of initial conditions based on incoming observations. Notably, the computational time of PINNYYF in each time step is substantially shorter than the observation sampling period.

PINNYYF was evaluated across three nonlinear problem cases, demonstrating high accuracy, computational efficiency, and reduced storage demands, suggesting its potential applicability to real-world nonlinear filtering tasks.

Despite the promising results, the applicability of PINNYYF to more complex, real-world scenarios remain challenging. High-dimensional problems with multiple state variables introduce significant difficulties for the FKE solver. Furthermore, improvements in both computational efficiency and accuracy are necessary for large-scale nonlinear filtering tasks.

Several potential directions for future research arise from this work. Firstly, extending the PINNYYF to higher-dimensional time-variant filtering problems is a priority. Secondly, applying the algorithm to real-world engineering problems using empirical data would provide further validation. Thirdly, exploring alternative neural network architectures, such as trunk-branch structures, could offer improvements in the precision of the approximate FKE solver.

\section*{Acknowledgements}
We appreciate the useful assistances and discussions from our colleague Dr. Xiaopei Jiao. This work is partially supported by Beijing Institute of Mathematical Sciences and Applications.

\bibliographystyle{utphys}
\bibliography{reference.bib}

\providecommand{\href}[2]{#2}\begingroup\raggedright\begin{thebibliography}{10}

\bibitem{Chen2009}
H.~Chen and K.~Chang.  ~Novel nonlinear filtering and prediction method for maneuvering target tracking. {\em IEEE Transactions on Aerospace and Electronic Systems} {\bfseries 45} no.~1~(2009), 237--249.

\bibitem{Date2011}
P.~Date and K.~Ponomareva.  ~Linear and non-linear filtering in mathematical finance: a review. {\em IMA Journal of Management Mathematics} {\bfseries 22} no.~3~(2011), 195--211.

\bibitem{Zhang2011}
Z.~Zhang, X.~Liu, G.~Yang,  and J.~Min.  ~Application and research on extended kalman prediction algorithm in target tracking system. in {\em Electrical, Information Engineering and Mechatronics 2011: Proceedings of the 2011 International Conference on Electrical, Information Engineering and Mechatronics (EIEM 2011)} pp.~1167--1173.
\newblock Springer.

\bibitem{Crisan2011}
D.~Crisan and B.~Rozovskii.  {\em The Oxford handbook of nonlinear filtering}.
\newblock Oxford University Press 2011.

\bibitem{Kutschireiter2020}
A.~Kutschireiter, S.~C. Surace,  and J.-P. Pfister.  ~The hitchhiker’s guide to nonlinear filtering. {\em Journal of Mathematical Psychology} {\bfseries 94} (2020), 102307.

\bibitem{Kalman1960}
R.~E. Kalman.  ~A new approach to linear filtering and prediction problems. {\em Journal of Basic Engineering} {\bfseries 82} no.~1~(1960), 35--45.

\bibitem{Kalman1961}
R.~E. Kalman and R.~S. Bucy.  ~New results in linear filtering and prediction theory. {\em Journal of Basic Engineering} {\bfseries 83} no.~1~(1961), 95--108.

\bibitem{Jazwinski2007}
A.~H. Jazwinski.  {\em Stochastic processes and filtering theory}.
\newblock Courier Corporation 2007.

\bibitem{Julier1997}
S.~J. {Julier} and J.~K. {Uhlmann}.  ~New extension of the kalman filter to nonlinear systems. in {\em Signal Processing, Sensor Fusion, and Target Recognition VI}, I.~{Kadar}, ed. vol.~3068 of {\em Society of Photo-Optical Instrumentation Engineers (SPIE) Conference Series} pp.~182--193.
\newblock 1997.

\bibitem{Evensen2003}
G.~Evensen.  ~The ensemble kalman filter: Theoretical formulation and practical implementation. {\em Ocean dynamics} {\bfseries 53} (2003), 343--367.

\bibitem{Gordon1993}
N.~Gordon, D.~Salmond,  and A.~Smith.  ~Novel approach to nonlinear/non-gaussian bayesian state estimation. {\em IEE Proceedings F (Radar and Signal Processing)} {\bfseries 140} (1993), 107--113.

\bibitem{Duncan1967}
T.~E. Duncan.  {\em Probability densities for diffusion processes with applications to nonlinear filtering theory and detection theory}.
\newblock Stanford University 1967.

\bibitem{Mortensen1966}
R.~E. Mortensen.  {\em Optimal control of continuous-time stochastic systems}.
\newblock University of California, Berkeley 1966.

\bibitem{Zakai1969}
M.~Zakai.  ~On the optimal filtering of diffusion processes. {\em Zeitschrift für Wahrscheinlichkeitstheorie und verwandte Gebiete} {\bfseries 11} no.~3~(1969), 230--243.

\bibitem{Yau2000}
S.-T. Yau and S.~S.-T. Yau.  ~Real time solution of nonlinear filtering problem without memory i. {\em Mathematical Research Letters} {\bfseries 7} no.~6~(2000), 671--693.

\bibitem{Yau2008}
S.-T. Yau and S.~S.-T. Yau.  ~Real time solution of the nonlinear filtering problem without memory ii. {\em SIAM Journal on Control and Optimization} {\bfseries 47} no.~1~(2008), 163--195.

\bibitem{Luo2013}
X.~Luo and S.~S.-T. Yau.  ~Complete real time solution of the general nonlinear filtering problem without memory. {\em IEEE Transactions on Automatic Control} {\bfseries 58} no.~10~(2013), 2563--2578.

\bibitem{Dong2020}
W.~Dong, X.~Luo,  and S.~S.-T. Yau.  ~Solving nonlinear filtering problems in real time by legendre galerkin spectral method. {\em IEEE Transactions on Automatic Control} {\bfseries 66} no.~4~(2020), 1559--1572.

\bibitem{Shi2018}
J.~Shi, Z.~Yang,  and S.~S.-T. Yau.  ~Direct method for yau filtering system with nonlinear observations. {\em International Journal of Control} {\bfseries 91} (2018), 678--687.

\bibitem{Yueh2014}
M.-H. Yueh, W.-W. Lin,  and S.-T. Yau.  ~An efficient algorithm of yau-yau method for solving nonlinear filtering problems. {\em Communications in Information and Systems} {\bfseries 14} no.~2~(2014), 111--134.

\bibitem{Wang2019}
Z.~Wang, X.~Luo, S.~S.-T. Yau,  and Z.~Zhang.  ~Proper orthogonal decomposition method to nonlinear filtering problems in medium-high dimension. {\em IEEE Transactions on Automatic Control} {\bfseries 65} no.~4~(2019), 1613--1624.

\bibitem{Shi2024}
J.~Shi, X.~Jiao,  and S.~S.-T. Yau.  ~Dglg: A novel deep generalized legendre-galerkin approach to optimal filtering problem. {\em IEEE Transactions on Automatic Control} (2024), 2584--2590.

\bibitem{Chen2024}
X.~Chen, Z.~Sun, Y.~Tao,  and S.~S.-T. Yau.  ~A uniform framework of yau-yau algorithm based on deep learning with the capability of overcoming the curse of dimensionality. {\em IEEE Transactions on Automatic Control} (2024), 339--354.

\bibitem{Tao2023}
Y.~Tao, J.~Kang,  and S.~S.-T. Yau.  ~Neural projection filter: Learning unknown dynamics driven by noisy observations. {\em IEEE Transactions on Neural Networks and Learning Systems} (2023), 9508--9522.

\bibitem{Raissi2019}
M.~Raissi, P.~Perdikaris,  and G.~E. Karniadakis.  ~Physics-informed neural networks: A deep learning framework for solving forward and inverse problems involving nonlinear partial differential equations. {\em Journal of Computational physics} {\bfseries 378} (2019), 686--707.

\bibitem{Mishra2023}
S.~Mishra and R.~Molinaro.  ~Estimates on the generalization error of physics-informed neural networks for approximating pdes. {\em IMA Journal of Numerical Analysis} {\bfseries 43} no.~1~(2023), 1--43.

\bibitem{Chen1995}
T.~Chen and H.~Chen.  ~Universal approximation to nonlinear operators by neural networks with arbitrary activation functions and its application to dynamical systems. {\em IEEE transactions on neural networks} {\bfseries 6} no.~4~(1995), 911--917.

\bibitem{Speidel2021}
J.~Speidel and J.~Speidel.  ~Randomly changing time-variant systems. {\em Introduction to Digital Communications} (2021), 165--175.

\bibitem{Pintelon2016}
R.~Pintelon, E.~Louarroudi,  and J.~Lataire.  ~Time-variant frequency response function measurement of multivariate time-variant systems operating in feedback. {\em IEEE Transactions on Instrumentation and Measurement} {\bfseries 66} no.~1~(2016), 177--190.

\bibitem{Olaleye2004}
F.~Olaleye, B.~Huang,  and E.~Tamayo.  ~Industrial applications of a feedback controller performance assessment of time-variant processes. {\em Industrial and engineering chemistry research} {\bfseries 43} no.~2~(2004), 597--607.

\bibitem{Abdi2010}
H.~Abdi and L.~J. Williams.  ~Principal component analysis. {\em Wiley interdisciplinary reviews: computational statistics} {\bfseries 2} no.~4~(2010), 433--459.

\bibitem{Davis1980}
M.~H. Davis and S.~I. Marcus.  ~An introduction to nonlinear filtering. in {\em Stochastic Systems: The Mathematics of Filtering and Identification and Applications: Proceedings of the NATO Advanced Study Institute held at Les Arcs, Savoie, France, June 22–July 5, 1980} pp.~53--75.
\newblock Springer.

\bibitem{Paszke2017}
A.~Paszke, S.~Gross, S.~Chintala, G.~Chanan, E.~Yang, Z.~DeVito, Z.~Lin, A.~Desmaison, L.~Antiga,  and A.~Lerer.  ~Automatic differentiation in pytorch.

\bibitem{Chen2017}
X.~Chen, X.~Luo,  and S.~S.-T. Yau.  ~Direct method for time-varying nonlinear filtering problems. {\em IEEE Transactions on Aerospace and Electronic Systems} {\bfseries 53} (2017), 630--639.

\bibitem{Kingma2014}
D.~P. Kingma and J.~Ba.  ~Adam: A method for stochastic optimization. 2017.
\newblock \url{https://arxiv.org/abs/1412.6980}.

\bibitem{He2015}
K.~He, X.~Zhang, S.~Ren,  and J.~Sun.  ~Deep residual learning for image recognition. in {\em Proceedings of the IEEE conference on computer vision and pattern recognition} pp.~770--778.

\end{thebibliography}\endgroup

\end{document}